\documentclass[11pt,reqno,a4paper]{amsart}
\usepackage{color}
\usepackage{amssymb,amsmath}
\usepackage{mathtools}
\usepackage[latin1]{inputenc}
\usepackage[active]{srcltx}
\usepackage{graphicx}
\usepackage{exscale,relsize}
\usepackage{textgreek}
\usepackage{epsfig,graphics}
\usepackage{psfrag}
\usepackage{caption}
\usepackage{subcaption}
\usepackage[toc,page]{appendix}
\usepackage{bigints}
\def\R{\mathbb{R}}
\def\m1{{I\!\!M}}

\newcommand{\Abracket}[1]{\left<#1\right>} 
\newcommand{\parenthesis}[1]{\left(#1\right)} 
\newcommand{\braces}[1]{\left\{#1\right\}} 

\renewcommand{\to}{\rightarrow}
\newcommand{\pa}{\partial}

\newcommand{\ino}{\int_{\Omega}}

\newcommand{\rife}[1]{(\ref{#1})}
\newcommand{\ov}[1]{\overline{#1}}

\newcommand{\sscp}{\scriptscriptstyle}

\renewcommand{\dfrac}{\displaystyle\frac}
\newcommand{\finedim}{\hspace{\fill}$\square$}
\newcommand{\intbar}{\mathop{\int\makebox(-15.5,0){\rule[6pt]{.7em}{0.3pt}}\kern-6pt}\nolimits}

\renewcommand{\i}{\infty}
\newcommand{\ii}{\infty}

\newcommand{\al}{\alpha}

\newcommand{\sg}{\sigma}

\newcommand{\om}{\Omega}
\newcommand{\lm}{\lambda}

\newcommand{\rl}{g_{\,\!\sscp \lm}}
\newcommand{\rla}{{g}_{\,\!\sscp \lm,\sscp \al}}

\newcommand{\rlq}{g'_{\sscp \lm}}

\newcommand{\rlqa}{g'_{\sscp \lm,\al}}

\newcommand{\thl}{\theta_{\sscp \lm}}

\renewcommand{\rho}{\mbox{\Large \textrho}}

\newcommand{\pl}{\psi_{\sscp \lm}}
\newcommand{\xil}{\ul}
\newcommand{\ul}{u_{\sscp \lm}}

\newcommand{\wl}{w_{\sscp \lm}}

\newcommand{\vl}{\eta_{\sscp \lm}}

\newcommand{\ssl}{\sscp \lm}
\newcommand{\ml}{m_{\sscp \lm}}
\newcommand{\all}{\al_{\ssl}}

\newcommand{\el}{E_{\ssl}}

\DeclareMathOperator{\Eigen}{Eigen}
\DeclareMathOperator{\Spect}{Spect}
\DeclareMathOperator{\Span}{Span}

\newtheorem{theorem}{Theorem}[section]
\newtheorem{proposition}[theorem]{Proposition}
\newtheorem{lemma}[theorem]{Lemma}
\newtheorem{corollary}[theorem]{Corollary}
\newtheorem{remark}[theorem]{Remark}
\newtheorem{definition}[theorem]{Definition}
\newcommand{\brm}{\begin{remark}\rm}
\newcommand{\erm}{\end{remark}}
\newcommand{\bdf}{\begin{definition}\rm}
\newcommand{\edf}{\end{definition}}
\newcommand{\bte}{\begin{theorem}}
\newcommand{\ete}{\end{theorem}}
\newcommand{\bpr}{\begin{proposition}}
\newcommand{\epr}{\end{proposition}}
\newcommand{\ble}{\begin{lemma}}
\newcommand{\ele}{\end{lemma}}
\newcommand{\bco}{\begin{corollary}}
\newcommand{\eco}{\end{corollary}}
\newcommand{\beq}{\begin{equation}}
\newcommand{\eeq}{\end{equation}}
\newcommand{\bdm}{\begin{displaymath}}
\newcommand{\edm}{\end{displaymath}}

\newcommand{\graf}[1]{\left\{\begin{array}{ll}#1\end{array}\right.}

\def\sideremark#1{\ifvmode\leavevmode\fi\vadjust{\vbox to0pt{\vss
 \hbox to 0pt{\hskip\hsize\hskip1em \vbox{\hsize2.1cm\tiny\raggedright\pretolerance10000 \noindent #1\hfill}\hss}\vbox to15pt{\vfil}\vss}}}

\begin{document}
\numberwithin{equation}{section}
\parindent=0pt
\hfuzz=2pt
\frenchspacing

\title[Bifurcation diagaram for Grad-Shafranov equation]{
Qualitative bifurcation diagram for \\ Grad-Shafranov type equations}


\author[D. Bartolucci]{Daniele Bartolucci}
\address{Daniele Bartolucci, Department of Mathematics, University of Rome \emph{"Tor Vergata"}, Via della ricerca scientifica n.1, 00133 Roma.}
\email{bartoluc@mat.uniroma2.it}

\author[G. Gu]{Guangze Gu}
\address{Guangze Gu, Department  of  Mathematics, Yunnan  Normal  University, Kunming 650500, China.}
\email{guangzegu@163.com}

\author[A. Jevnikar]{Aleks Jevnikar}
\address{Aleks Jevnikar, Department of Mathematics, Computer Science and Physics, University of Udine, Via delle Scienze 206, 33100 Udine, Italy.}
\email{aleks.jevnikar@uniud.it}

\author[J. Wei]{Juncheng Wei}
\address{Juncheng Wei, Department of Mathematics, Chinese University of Hong Kong, Shatin, Hong Kong.}
\email{wei@math.cuhk.edu.hk}

\author[R. Wu]{Ruijun Wu}
\address{Ruijun Wu, School of mathematics and statistics, Beijing Institute of Technology, Zhongguancun South Street No. 5, 100081 Beijing, P.R.China.}
\email{ruijun.wu@bit.edu.cn}

\thanks{D.B.
is partially supported by the MIUR Excellence Department Project MatMod@TOV
awarded to the Department of Mathematics, University of Rome ``Tor Vergata'' and by the E.P.G.P. Project sponsored by the University of Rome ``Tor Vergata''. \\
\indent D.B. and A.J. are partially supported by INdAM-GNAMPA project ``{\em Aspetti qualitativi e analisi di blow-up per problemi differenziali ellittici}''.\\
\indent G.G. is supported by Xingdian talent support program of Yunnan Province, National Natural Science Foundation of China (12261107, 12561020) and Yunnan Fundamental Research Projects (202601AT070048).\\
\indent J.W. is partially supported by GRF fund of RGC of Hong Kong entitled
``\emph{New frontiers in singularity formations of nonlinear partial differential equations}.''
}

\begin{abstract}
We study the qualitative behavior of solutions of
Grad-Shafranov type equations arising in plasma physics with general differential
operators and general nonlinearities. In particular, we extend recent estimates about
threshold values for uniqueness, monotonicity and non-existence of the free boundary. The argument is based on a refined spectral analysis for weighted non-local problems together with comparison techniques and level set analysis.
\end{abstract}
\maketitle

\

\section{\bf Introduction}

\medskip

Let $\om\subset\R^N$ be a smooth bounded domain with $N\geq2$ and let $\lm\geq0$ be fixed. We study here solutions $(\al,\psi)\in\R\times C^{2,r}(\ov \om)$, $r\in(0,1)$, to the following Grad-Shafranov type problem
\beq \label{prob}
\begin{cases}
\mathcal D \psi = g(x,\al+\lm\psi) \quad \mbox{in } \om, \\ \\
\bigintss_\om g(x,\al+\lm\psi)= 1, \\ \\
\psi=0 \quad \mbox{on } \partial\om.
\end{cases}
\eeq
Here, $\mathcal D$ is the operator defined by
$$
	\mathcal D\psi= -\sum_{i,j=1}^N \dfrac{\pa}{\pa x_i}\left( a_{ij}(x)\dfrac{\pa \psi}{\pa x_j} \right)
$$
with $a_{ij}\in C^2(\ov\om)$, $a_{ij}=a_{ji}$, which is assumed to be uniformly elliptic in $\ov\om$, that is
\beq\label{ell}
	\sum_{i,j=1}^N a_{ij}(x) \xi_i \xi_j \geq A|\xi|^2 \qquad \forall x\in\ov\om, \ \forall\xi\in\R^N,
\eeq
for some $A>0$.

\medskip

The continuous function $g=g(x,z):\ov\om\times\R\to[0,+\i)$, $C^1$ in $z$, has the usual properties
\begin{align*}
&g(\cdot,z)=0 \quad \mbox{for } z\leq0, \\
&g_z(\cdot,z)>0 \quad \mbox{for } z>0.
\end{align*}
Moreover, we assume the following, for $z>0$:
\begin{equation}
\dfrac{g_z(\cdot,z)}{g(\cdot,z)}>\dfrac1z, \label{g1}
\end{equation}

\smallskip

\begin{equation}  \label{g2}
\dfrac{g_z(\cdot,z)}{g(\cdot,z)}\leq \dfrac pz \quad \mbox{and} \quad g(\cdot,z)\leq z^p,
\end{equation}
for some $p\in(1,p_N)$, where
$$
p_N= \begin{cases}
+\i, &\quad N=2, \\
\dfrac{N}{N-2}, &\quad N\geq3.
\end{cases}
$$
For the sake of clarity, and for future references, the model case is given by
\beq \label{model}
\mathcal D=-\Delta, \quad g(x,z)=[z]_+^p,
\eeq
which has been very much studied in the literature. Here, $A=1$. With respect to \eqref{model}, the assumptions \eqref{g1} and \eqref{g2} are encoding the superlinear ($p>1$) and subcritical ($p<p_N$) behavior of the nonlinearity, respectively. The condition \eqref{g2} can be replaced by $g(\cdot,z)\leq Cz^p$ for some $C>0$, but we stick with the present formulation for the sake of simplicity.

\

System \eqref{prob} is related to Grad-Shafranov type equations arising in the physics of Tokamak's plasma, see \cite{Frei,Kad,Stac} and the discussion in the appendix. After the first pioneering works for the linear case \cite{Te,Te2} and the general subcritical case \cite{BeBr}, there has been a lot of effort to study existence, uniqueness, multiplicity of solutions, structure of the free boundary, mainly for the model problem \eqref{model}, and we refer the interested readers to \cite{BJW2} for a detailed list of references. In particular, the threshold $p_N$ turns out to be a critical growth for this class of problems with respect both to a priori estimates and existence of solutions, see for example \cite{BeBr,Ort}.

\medskip

On the other hand, the qualitative study of the shape of branches of solutions has been initiated just recently \cite{BHJY, BJ1, BJWW1, BJWW2,BJWW3,BJW2} for the model case \eqref{model}, see also \cite{L-E} for a related problem. In particular, relevant advances has been recently made in \cite{BJWW1,BJW2},
where one can find sharp estimates about threshold values for uniqueness, monotonicity
and non-existence of the free boundary. However, in the physical applications the exact form of the nonlinearity $g$ is an unknown of the problem and we refer to \cite{bbf} and the references quoted therein for this aspect, called ``reconstruction problem''. This motivates the rigorous analysis of \eqref{prob} under very general assumptions both on the differential operator and on $g$. It is our aim here to make a first step in this direction.

\

Before stating our results we make some preparations. We will write $(\all, \pl)$ for solutions of \eqref{prob} when we want to stress the dependence on the parameter $\lm$. For simplicity we assume $|\om|=1$ and for fixed $t\geq 1$ we denote
\beq \label{sob}
\Lambda(\om,t)=\inf\limits_{w\in H^1_0(\om), w\equiv \!\!\!\!/ \;0}\dfrac{\ino |\nabla w|^2}{\left(\ino |w|^{t}\right)^{\frac2t}}\,,
\eeq
which provides the best constant in the Sobolev embedding
$$
\|w\|_p\leq S_p(\om)\|\nabla w\|_2, \qquad S_p(\om)=\Lambda^{-1/2}(\om,p), \qquad p\in[1,2p_{N}).
$$

\medskip
The energy of a solution $\psi$ of \eqref{prob} will be defined as
\beq\label{energy}
\el= \frac12  \ino \sum_{i,j=1}^N a_{ij}\, D_i\pl D_j \pl .
\eeq

\medskip

Clearly, for $\lm=0$ \eqref{prob} admits a unique solution. Our first task is to show that it is possible to extend this branch to a full explicit range of the parameter $\lm$, describing also qualitatively the shape of the diagram. This is done via bifurcation analysis and we thus associate to a solution $(\all,\pl)$ the following linearized equation
\beq\label{linear}
L_{\ssl}[\phi]=\mathcal D\phi-\lm \rlq [\phi]_{\ssl},
\eeq
where $g'_\lm=g_z(x,\all+\lm\pl)$ and
$$
<\phi>_{\ssl}=\dfrac{\ino g'_\lm  \phi}{\ino g'_\lm}, \quad \quad  [\phi]_{\ssl}=\phi \,-<\phi>_{\ssl},
$$
for any $\phi\in L^2(\Omega)$. We thus introduce the following definition.

\medskip

{\bf Definition.}
{\it $\sg=\sg(\all,\pl)\in\R$ is an eigenvalue of $L_{\ssl}$ if
\beq\label{lineq0}
\mathcal D \phi-\lm \rlq [\phi]_{\ssl}=\sg\rlq [\phi]_{\ssl},
\eeq
has a non-trivial weak solution $\phi\in H^1_0(\om)$.}

\medskip

This is the natural spectral setting and we postpone the discussion and all the details to section \ref{sec:spectral}. Here we just point out that the operator \eqref{linear} involves both nonlocal terms $[\phi]_{\ssl}$ and a weight $\rlq$ which may vanish on a large set. Therefore, some care is needed to build the spectral framework. Nevertheless, we can show that the first associated eigenvalue $\sg_1(\all,\pl)$ is well-defined and introduce the threshold
$$
\lm_*(\om,p)=\sup\bigr\{\lm>0\,:\,\sg_1(\alpha_{\mu},\psi_{\mu})>0\,\mbox{\rm \,for any solution of } \eqref{prob}_{ \mu},\,\forall\,\mu<\lm\bigr\}.
$$
Then, we prove the following. Here, $\mathcal D$ and $g$ satisfy the assumptions after \eqref{prob} and we recall also the constant $A$ in \eqref{ell}.
\bte\label{thm2}
Let $N\geq2$ and $p\in (1,p_{N})$. Then, it holds:

\medskip

\begin{itemize}
\item[1.] \emph{(Spectral estimate):} $\lm_*(\om,p)>\frac{A}{p} \Lambda(\om,2p)$.

\

\item[2.] \emph{(Uniqueness):} For any $\lm<\lm_*(\om,p)$ there exists a unique solution $(\all, \pl)$ of \eqref{prob}.

\

\item[3.] \emph{(Qualitative behavior):} For any $\lm<\lm_*(\om,p)$ we have
$$
\frac{d \el}{d\lm}>0.
$$
Moreover, for any $\lm\leq\frac{A}{p} \Lambda(\om,2p)$ we have
$$
\frac{d \all}{d\lm}<0.
$$
\end{itemize}
\ete

\medskip

More precisely, we get a unique $C^1$ branch of solutions, starting from $\lm=0$, for which $\sg_1(\all,\pl)>0$, where the pointwise monotonicity of $\pl$ is replaced by the monotonicity of $\el$ and $\all$. The different range for which we get monotonicity of $\all$ is due to the fact that we are not able, up to now, to fully exploit the new modified spectral setting and we have to rely on more standard spectral estimates, see Proposition \ref{pr3.2.best}. It is a challenging open problem to study whether the monotonicity of $\all$ holds up to the threshold $\lm_*(\om,p)$. Remark that this is true for the model problem as recently shown in \cite{BJWW1}.

\

We next turn to the study of the so called free boundary $\partial\{x\in\Omega \,:\, \alpha+\lm\psi>0\}$. Since by the maximum principle $\psi\geq0$ in $\Omega$, the existence of the free boundary depends on the sign of $\alpha$. This motivates  the following definition.

\medskip

{\bf Definition.} \emph{We say that a solution to \eqref{prob} is non-negative (resp. positive) if $\al\geq0$ (resp. $\al>0$). }

\medskip

For $\lambda=0$ we have a unique solution which is also positive, and then we deduce from Theorem \ref{thm2} that the following quantity, which we call positivity threshold, is well-defined and strictly positive,
\beq \label{positive}
\lm_+(\om,p)=\sup\{\lm>0\,:\,\al_{\mu}>0\,\mbox{\rm \,for any solution of } \eqref{prob}_{ \mu},\,\forall\,\mu<\lm\}.
\eeq
Here and in the sequel we are still assuming $\mathcal D$ and $g$ satisfy the assumptions after \eqref{prob}. Then, we have the following.
\bte\label{thm1}
Let $N\geq2$ and $p\in (1,p_{N})$. Then, it holds
$$
\lm_+(\om,p)> \frac{A}{p} \Lambda(\om,2p).
$$
\ete

\medskip

In particular, for any $\lambda\leq\frac{A}{p} \Lambda(\om,2p)$ we get a unique branch of solutions which are positive, i.e. with no free boundary.

\medskip

For $N=2$ we can improve the estimate in Theorem \ref{thm1} assuming that we have a non-degenerate bound from below on the $x$ variable:
\beq \label{extra}
C_1\,\tilde g(z) \leq g(\cdot,z)\leq C_2\,\tilde g(z),
\eeq
for any $z$, for some constants $C_1,\,C_2>0$ and $\tilde g(z)$ satisfying \eqref{g1}-\eqref{g2}. For example, one may consider the model
$$
g(x,z)=f(x)\tilde g(z), \quad C_1 \leq f(x)\leq C_2.
$$
We prove the following.
\bte\label{thm3}
Let $N=2$ and $p\in (1,p_{N})$. Assume also \eqref{extra}. Then, it holds
$$
\lm_+(\om,p)> A\left(\frac{C_1}{C_2}\right)^{\frac{p-1}{2p}}\left(\dfrac{8\pi}{p+1}\right)^{\frac{p-1}{2p}} \Lambda^{\frac{p+1}{2p}}(\om,p+1).
$$
\ete

\medskip

We conjecture that
$$
\left(\frac{8\pi}{p+1}\right)^{\frac{p-1}{2p}}\Lambda^{\frac{p+1}{2p}}(\om,p+1)>\frac{1}{p} \Lambda(\om,2p).
$$
This was confirmed in \cite{BJW2} for $p$ large.

\

These are the first results about the qualitative behavior of the branch of solutions for general differential operators and general nonlinearities and they hold for any smooth and bounded domain, in any dimension.
They generalizes the findings in \cite{BJ1, BJWW1, BJW2} related to the model case \eqref{model}. Let us also mention that in the latter simpler case, the result extends to $p=1$ for which the positivity threshold in Theorems \ref{thm1}, \ref{thm3} (here $A=1$ and $C_1=C_2=1$) $\Lambda(\om,2)$ is sharp, see for example \cite{BeBr,Pudam} and the well-known result by Temam \cite{Te2}.

\

We refine here the analysis started in \cite{BJ1, BJWW1} in order to handle the general problem \eqref{prob} which presents several subtle and technical point. The argument is based on bifurcation analysis, the core being the spectral study of the linearized equation. To catch the effect of the constraint in \eqref{prob} we need to introduce a non-standard spectral set up. This affects the classical results, for instance the maximum principle fails and, moreover, the first eigenvalue is not simple in general. We refer to \cite{BJ0} for some explicit examples concerning these aspects in a similar setting. Furthermore, recall that the weight may vanish on a large set. A detailed study of how this affects the spectral properties is thus needed. Nevertheless, it turns out that it nicely describes the branch of solutions outside the kernel of the linearized equation. This yields in particular the uniqueness of solutions. The shape of the bifurcation diagram is then based on the spectral decomposition along the branch. The threshold in Theorem \ref{thm1} is obtained via some comparison argument. Finally, the improved estimate in Theorem \ref{thm3} is deduced from a universal energy estimate, of independent interest, which in turn is based on a level set analysis.

\

The article is organized as follows. In section \ref{sec:spectral} we introduce the spectral set up and prove some related key results, in section \ref{sec:bif} we prove Theorem \ref{thm2}, i.e. the spectral estimate together with uniqueness and monotonicity properties and finally, section \ref{sec:positive} is devoted to the proof of the estimates concerning the positivity threshold, namely Theorems \ref{thm1}, \ref{thm3}. Some considerations about Grad-Shafranov type equations are collected in the appendix.

\
	
\section{\bf Spectral set up} \label{sec:spectral}

\medskip

In this section we introduce the spectral setting and prove its basic properties needed to study \eqref{prob}. When there is no risk of confusion, for a solution $(\all,\pl)$ of \eqref{prob}, with a little abuse of notation, we simply denote
$$
	g_{\lm,\al} =g_{\lm,\al}(\psi)=g(x,\al+\lm\psi), \quad g_{\lm}=g(x,\all+\lm\pl)
$$
and
$$
	 g'_{\lm,\al}= g_z(x,\al+\lm\psi), \quad g'_\lm=g_z(x,\all+\lm\pl).
$$
We will also write
$$
	\ml=\ino g'_{\lm}.
$$
Moreover, for $\mu,\varphi\in L^2(\om)$, we define
$$
<\mu>_{\ssl}=\dfrac{\ino g'_\lm  \mu}{\ino g'_\lm}, \quad \quad  [\mu]_{\ssl}=\mu \,-<\mu>_{\ssl},
$$
and
$$
<\mu,\varphi>_{\ssl}=\dfrac{\ino g'_\lm  \mu\varphi}{\ino g'_\lm}.
$$

\medskip

Let $C^{2,r}_0(\ov{\om}\,)$ be the space of $C^{2,r}(\ov{\om}\,)$ functions with zero boundary conditions. We then define
$$
 F\colon (-1,+\infty)\times \R \times C^{2,r}_0(\ov{\Omega}) \to C^r(\ov{\Omega}),
$$
\begin{align}
&F(\lm,\al,\psi)=\mathcal D \psi -g_{\lm,\al}(\psi) \label{eF}
\end{align}
and
$$
\Phi:(-1,+\infty)\times \R \times C^{2,r}_0(\ov{\Omega})\to  \R\times C^{r}(\ov{\om}\,)
$$

$$
\Phi(\lm,\al,\psi)=\left(\begin{array}{cl}-1+\ino g_{\lm,\al}\\ \\F(\lm,\al,\psi)\end{array}\right).
$$
Now,
we consider
$$
D_{\al,\psi}\Phi(\lm,\al,\psi):\R \times C^{2,r}_0(\ov{\om}\,) \to  \R\times C^{r}(\ov{\om}\,),
$$

$$
D_{\al,\psi}\Phi(\lm,\al,\psi)[s,\phi]=
\left(\begin{array}{cr}
\ino \left(D_\psi \rla[\phi]+d_{\al}\rla [s]\right)\\ \\
D_\psi F(\lm,\al,\psi)[\phi]+d_{\al}F(\lm,\al,\psi)[s] \end{array}\right),
$$
where, for $\phi \in C^{2,r}_0(\ov{\om}\,)$,
$$
D_\psi F(\lm,\al,\psi)[\phi]=
\mathcal D  \phi -\lm \rlqa \phi,
$$

$$
D_\psi \rla[\phi]=\lm\rlqa \phi,
$$
and, for $s \in \R$,

$$
d_{\al}F(\lm,\al,\psi) [s]=- \rlqa s,
$$

$$
d_{\al}\rla [s]= \rlqa s.
$$

\

Clearly, solutions $(\all,\pl)$ to \eqref{prob} correspond to
$\Phi(\lm,\all,\pl)=(0,0)$ and we consider the linearized equation
$$
L_{\ssl}[\phi]=\mathcal D\phi-\lm \rlq [\phi]_{\ssl}
$$
and we recall that $\sg=\sg(\all,\pl)\in\R$ is said to be an eigenvalue of $L_{\ssl}$ if
$$
\mathcal D \phi-\lm \rlq [\phi]_{\ssl}=\sg\rlq [\phi]_{\ssl},
$$
has a non-trivial weak solution $\phi\in H^1_0(\om)$. The corresponding eigenspace is denoted by $\Eigen(L_{\ssl};\sg)$. We introduce hereafter the spectral setting taking into account the nonlocal terms and the weight which may vanish on a large set, recalling some classical facts for weighted operators from \cite{MaMi}.

\medskip

We rewrite the eigenvalue equation as
\begin{align}\label{eq:weighted eigenvalue}
    \mathcal D \phi = (\lm+\sigma)\rlq [\phi]_{\ssl}
\end{align}
and use the Green's function $G$ for~$\mathcal D$ with Dirichlet boundary condition to get
$$
    \phi =  \frac{\lm+\sg}{\lm}\;  G* \parenthesis{\lm\rlq[\phi]_{\ssl}}.
$$
We now consider the operator~$T_{\ssl} \colon H^1_0(\Omega)\to H^1_0(\Omega)$
$$
    T_{\ssl}(\phi) \coloneqq G* \parenthesis{\lm\rlq[\phi]_{\ssl}}
$$
to rewrite the eigenvalue equation as
\begin{align}\label{eq:eigenvalue for T}
    T_{\ssl}(\phi)= \mu\phi,
\end{align}
where
$$
\mu = \frac{\lm}{\lm+\sg}.
$$
Therefore, there is a one to one correspondence between the spaces $\Eigen(L_{\ssl};\sigma)$ and $\Eigen(T_{\ssl};\frac{\lm}{\lm+\sg})$.
We focus now on~$T_{\ssl}$. It is easy to check that it is a linear self-adjoint compact operator on ~$H^1_0(\Omega)$ equipped with the associated inner product
$$
    \Abracket{\xi,\eta}_{H^1_0}:=\int_{\Omega} \sum_{i,j}a_{ij}D_i\xi D_j\eta, \qquad  \forall \xi,\eta\in H^1_0(\Omega).
$$
For instance,
\begin{align*}
    \Abracket{T_{\ssl}(\xi),\; \eta}_{H^1_0}
    =& \int_{\Omega} \sum_{i,j}a_{ij}D_i(T_{\ssl}(\xi)) D_j\eta
    =\int_{\Omega} (\mathcal D T_{\ssl}(\xi))\eta\\
    =&\int_{\Omega} \lm\rlq [\xi]_{\ssl} \eta
    =\int_{\Omega} \lm\rlq [\xi]_{\ssl} [\eta]_{\ssl} \\
		=&\Abracket{\xi,\; T_{\ssl}(\eta)}_{H^1_0}.
\end{align*}
Let us denote by $\Spect(T_{\ssl})$ the spectrum of~$T_{\ssl}$. It follows by standard results that it has countably many real nonzero eigenvalues with finite multiplicity plus possibly zero. Moreover, zero is the unique accumulation point.
Since~$\lm+\sigma\geq 0$, one has~$\Spect(T_{\ssl})\subset\R_+$. Let us write
$$
    \mu_1\geq \mu_2\geq \mu_3\geq \cdots >0,
$$
with $\mu_j \to 0$ as $j\to+\infty$ and denote the corresponding eigenfunctions by~$\phi_j$,~$j\in\mathbb{N}$, which satisfy
$$
    \delta_{jk}=\Abracket{\phi_j,\phi_k}_{H^1_0} \qquad \forall j,k\geq 1.
$$
Following \cite{MaMi} we consider the bilinear form,
$$
    \mathcal{B}(\xi,\eta)\coloneqq \lm\int_{\om} \rlq [\xi]_{\ssl}[\eta]_{\ssl}, \qquad \forall \xi,\eta\in H^1_0(\Omega),
$$
and we know that
\begin{align*}
                \mu_1= &\max \braces{ \mathcal{B}(\phi,\phi) \mid \phi\in H^1_0(\Omega),\; \Abracket{\phi,\phi}_{H^1_0}=1}  \\
                =&\; \mathcal{B}(\phi_1,\phi_1),
\end{align*}
\begin{align*}
                \mu_j=& \max\braces{\mathcal{B}(\phi,\phi)\mid \phi\in H^1_0(\Omega),\; \Abracket{\phi,\phi}_{H^1_0}=1, \; \Abracket{\phi,\phi_k}_{H^1_0} =0, \; 1\leq k\leq j-1} \\
                =& \;\mathcal{B}(\phi_j,\phi_j).
\end{align*}
Let
$$
\mathcal{H}_1\coloneqq \ov{\oplus_{j\geq 1} \Eigen(T_{\ssl};\mu_j) } = \ov{\Span\braces{\phi_j\mid j\geq 1}}.
$$
We collect in the following proposition the properties that we will use later on.
\bpr\label{pro-proj}
Let $\alpha<0$. Then,
$$
H^1_0(\Omega)= \Eigen(T_{\ssl};0)\oplus\mathcal{H}_1
$$
and there exist projections
$$
        P_0\colon H^1_0(\Omega)\to\Eigen(T_{\ssl};0),
$$
$$				
        P_1\colon H^1_0(\Omega)\to \mathcal{H}_1,
$$
such that for~$\psi\in H^1_0(\Omega)$,
    \begin{align} \label{decomposition}
        \psi= P_0 \psi + P_1 \psi = P_0\psi + \sum_{j=1}^{+\infty} \beta_j\phi_j
    \end{align}
    with
    \begin{align*}
        \beta_j
        =& \Abracket{\psi,\phi_j}_{H^1_0}
        = \int_{\Omega} \sum_{i,j}a_{ij}D_i\psi D_j\phi_j  \\
        =&\ino \frac{\lm}{\mu_j} \rlq [\phi_j]_{\ssl} \psi
        =\ino \frac{\lm}{\mu_j} \rlq [\phi_j]_{\ssl} [\psi]_{\ssl}.
    \end{align*}
Moreover, $0\in\Spect(T_{\ssl})$ if and only if~$\alpha<0$.
\epr

\proof

We start by proving the last property and we observe that by \eqref{eq:eigenvalue for T} we get
\begin{align*}
            \mu\mathcal D\phi= \lm\rlq[\phi]_{\ssl}.
\end{align*}
Denote by $\overline{\Omega_+}$ the support of $g_{\lm}$. For $\mu=0$ one has $\rlq[\phi]_{\ssl}\equiv 0$ a.e. in~$\Omega$ and hence $[\phi]_{\ssl}\equiv 0$ a.e. in~$\Omega_+$, i.e. $\phi = \Abracket{\phi}_{\ssl}$ is constant a.e. in~$\Omega_+$. Now, if~$\alpha\geq 0$, then~$\Omega_+= \Omega$ by the maximum principle and thus $\phi=0$ since~$\phi\in H^1_0(\Omega)$.

\medskip

On the other hand, if~$\alpha<0$ then~$\Omega\setminus \ov{\Omega_+}$ is a proper subset and any~$0\neq \phi \in H^1_0(\Omega\setminus \ov{\Omega_+})$ belongs to $\Eigen(T_{\ssl};0)$.

\medskip

By standard results we then have the orthogonal decomposition together with the associated projections $P_0$, $P_1$ with the above properties. The explicit expression of the Fourier coefficients follows by our construction.
\finedim

\

\section{\bf Bifurcation diagram} \label{sec:bif}

\medskip

In this section we collect all the results concerning the bifurcation diagram of \eqref{prob}, proving in particular uniqueness and monotonicity of its solutions and the spectral estimate. We consider here $\mathcal D$ and $g$ with the assumptions after \eqref{prob}.

\

\subsection{\bf Spectral estimate} \, \newline

Recalling the spectral setting of the previous section, we have
\begin{align*}
    \frac{1}{\mu_1}=\inf \left\{\frac{1}{\mathcal{B}(\phi,\phi)}\,:\,\phi\in H^1_0(\om),\; \Abracket{\phi,\phi}_{H^1_0}=1\right\},
\end{align*}
and then
\begin{align} \label{4.1}
\sg_1=\sg_1(\all,\pl)=\lm\left(\frac{1}{\mu_1}-1\right)=\inf\limits_{\phi \in H^1_0(\om)\setminus \{0\}}
\dfrac{\ino \sum_{i,j}a_{ij}D_i\phi D_j\phi - \lm \ino \rlq [\phi]_{\ssl}^2 }{\ino \rlq [\phi]_{\ssl}^2}.
\end{align}
On the other hand, the standard first eigenvalue is given by
\begin{align*}
     \nu_{1}(\all,\pl)\coloneqq \inf\limits_{w\in H^{1}_0(\om)\setminus\{0\}}
        \dfrac{\ino \sum_{i,j}a_{ij}D_i w D_j w-\lm \ino \rlq w^2  }{\ino \rlq w^2}.
\end{align*}
We have the following relation.
\ble \label{lem:xi1s} It holds
\beq\label{xi1s}
\sg_1(\all,\pl)>\nu_{1,\lm}.
\eeq
\ele

\proof
We have,
$$
\dfrac{\ino \sum_{i,j}a_{ij}D_i\phi D_j\phi - \lm \ino \rlq [\phi]_{\ssl}^2 }{\ino \rlq [\phi]_{\ssl}^2}=
\dfrac{\ino \sum_{i,j}a_{ij}D_i\phi D_j\phi }{\ino \rlq [\phi]_{\ssl}^2}-\lm=
$$

\smallskip

$$
\frac{1}{\ml}\dfrac{\ino \sum_{i,j}a_{ij}D_i\phi D_j\phi }{<[\phi]_{\ssl}^2>_{\ssl}}-\lm=\frac{1}{\ml}\dfrac{\ino \sum_{i,j}a_{ij}D_i\phi D_j\phi }{<\phi^2>_{\ssl}-<\phi>^2_{\ssl}}-\lm\geq
$$

\smallskip

$$
\frac{1}{\ml}\dfrac{\ino \sum_{i,j}a_{ij}D_i\phi D_j\phi }{<\phi^2>_{\ssl}}-\lm=\dfrac{\ino \sum_{i,j}a_{ij}D_i\phi D_j\phi}{\ino \rlq \phi^2}-\lm=
$$

\smallskip

$$
\dfrac{\ino \sum_{i,j}a_{ij}D_i\phi D_j\phi - \lm \ino \rlq \phi^2 }{\ino \rlq \phi^2}.
$$
The equality holds if and only if $<\phi>_{\ssl}=0$ and thus if and only if any
eigenfunction $\phi_1$ of $\nu_{1,\lm}$ satisfies $<\phi_1>_{\ssl}=0$, which is impossible since $\phi_1$ has constant sign.
\finedim

\

The following result provides then the spectral estimates of Theorem \ref{thm2}.
\bpr\label{preigen.LE}
If $(\all,\xil)$ is a solution of \eqref{prob} with $\lm\leq \frac{A}{p}\Lambda(\om,2p)$, we have
$$
\nu_{1,\lm}\geq0
$$
and then, in particular,
$$
\sg_1(\all,\pl)>0.
$$
\epr
\proof
Let $w\equiv \!\!\!\!\!/ \;0$. Then, we first use \eqref{ell} to handle the differential operator, then \eqref{g2} to estimate the nonlinearity $g$ and finally the constraint in \eqref{prob} to get
\begin{align*}
\dfrac{ \ino \sum_{i,j}a_{ij}D_i w D_j w }{\ino \rlq w^2} &\geq A \dfrac{\ino |\nabla w|^2}{\ino \rlq w ^2} \geq \frac{A}{p} \dfrac{\ino |\nabla w|^2}{\ino g_\lm^{({p-1})/{p}} w ^2} \\ & \\
	&\geq \frac{A}{p} \dfrac{\ino |\nabla w|^2}{\left( \ino g_\lm\right)^{\frac{p-1}{p}} \left(\ino w^{2p}\right)^{\frac1p}} = \frac{A}{p} \dfrac{\ino |\nabla w|^2}{\left(\ino w^{2p}\right)^{\frac1p}} \geq \frac{A}{p} {\Lambda(\om,2p)}\,.
\end{align*}
Therefore,
$$
\nu_{1,\lm}\geq \frac{A}{p}{\Lambda(\om,2p)}-\lm
$$
and the thesis follows.
\finedim

\

\subsection{\bf Uniqueness of solutions} \, \newline

We now discuss the uniqueness property in Theorem \ref{thm2}. Let us start with the following preliminary well-known estimate, see for example \cite{BeBr,BJ1}.
\ble\label{lemE1} Let $\ov{\lm}>0$ be fixed. Then, there
exists $C=C(r,\om,\ov{\lm},p,N)>0$  such that $\|\pl\|_{C^{2,r}_0(\ov{\om})}\leq C$ for any solution
$(\all,\pl)$ of $\eqref{prob}_{\lm}$ with $\lm\in [0,\ov{\lm}\,]$.
\ele

\proof
Since $\ino g(x,\all+\lm\pl)=1$, then for any $t\in [1,\frac{N}{N-1})$ there exists $C=C(t,N,\om)$ such that $\|\pl\|_{W_0^{1,t}(\om)}\leq C(t,N,\om)$ for any solution of \eqref{prob}, see \cite{St4}.
Then, the Sobolev inequalities give $\|\pl\|_{L^s(\om)}\leq C(s,N,\om)$, for any $1\leq s< \frac{N}{N-2}$,
for some $C(s,N,\om)$. It is also easy to see that $\all$ is bounded from above by looking at the constraint in \eqref{prob}. Now, by using assumption \eqref{g2} and recalling that $p<p_{N}$, we get the existence of some $m>1$ depending on $p$ and $N$ such that
$\|g(x,\all+\lm\pl)\|_{L^m(\om)}\leq C(p,N,\ov{\lm},s,\om)$,
for any $\lm\leq \ov{\lm}$,
for some $C(p,N,s,\ov{\lm},\om)$ and then, by standard elliptic estimates,
we conclude that $\|\pl\|_{W_0^{2,m}(\om)}\leq C_0(p,N,s,\ov{\lm},\om)$,
for any $\lm\leq \ov{\lm}$, for some
$C_0(p,N,s,\ov{\lm},\om)$.
At this point, since $\om$ is of class $C^{2,r}$, the conclusion follows by standard elliptic estimates
and a bootstrap argument.\finedim

\

The proof about uniqueness of solutions is based on the implicit function theorem and on the following result about the linearized problem and the branch of solutions to \eqref{prob} away from the spectrum.

\bpr\label{pro-lin}
Let $\lm_0\geq 0$ and let $(\al_{\sscp \lm_0},\psi_{\sscp \lm_0})$ be a solution of \eqref{prob} and suppose that $0\notin\mbox{\emph{Spect}}(L_{\lm_0})$. Then, locally around $(\lm_0,\al_{\sscp \lm_0},\psi_{\sscp \lm_0})$, the set of solutions of \eqref{prob} is a $C^1$ curve of solutions $\lm\mapsto (\all,\pl)$.

\medskip

In particular, for any $\lm<\lm_*(\om,p)$ there exists a unique solution $(\all, \pl)$ of \eqref{prob}.
\epr

\proof
The solutions of \eqref{prob} are related to the map $\phi$ introduced right after \eqref{eF}. We show that $D_{\al,\psi}\Phi(\lm_0, \al_{\sscp \lm_0}, \psi_{\sscp \lm_0})$ is an isomorphism. The conclusion then follows by the implicit function theorem. Let us write $\lm_0=\lm$.

\medskip

We first show Ker$(D_{\al,\psi}\Phi(\lm,\all,\pl))=0$. Indeed, suppose on the contrary that there exists a non-vanishing pair $(s,\phi)\in \R \times C^{2,r}_0(\ov{\om}\,) $ such that
$$D_{\al,\psi}\Phi(\lm,\all,\pl)[s,\phi]=(0,0).$$
Then the constraint $\left.\ino \left(D_\psi \rla[\phi]+d_{\al}\rla [s]\right)\right|_{(\all,\pl)}=0$ yields $s=- \lm<\phi>_{\ssl}$ and then we conclude that
$L_{\ssl}[\phi]=D_\psi F(\lm,\all,\pl)[\phi]+d_{\al}F(\lm,\all,\pl) [s_{\ssl}]=0$, i.e. $\phi$ is a non-trivial eigenfunction corresponding to the $0$ eigenvalue, a contradiction.

\medskip

We next prove that $D_{\al,\psi}\Phi(\lm,\all,\pl)$ is surjective and show that for any $(t,f)\in \R\times C^{r}(\ov{\om}\,)$,
$$
D_{\al,\psi}\Phi(\lm,\all,\pl)[s,\phi] =(t,f),
$$
has a solution $(s,\phi)\in \R\times C^{2,r}_0(\ov{\om}\,)$. Reasoning as above, from the first equation we get
\begin{align*}
     s = -\lm\Abracket{\phi}_{\ssl} + \frac{t}{\ino\rlq }
\end{align*}
and then the second gives
\begin{align}
    \mathcal D\phi- \lm \rlq [\phi]_{\ssl} = f+t\frac{\rlq}{\ino\rlq}.
\end{align}
and we want to prove that there is a unique~$\phi\in C^{2,r}_0(\ov{\Omega})$ solving this equation. This follows by showing that for $0\notin \Spect(L_{\ssl}) $, $L_{\ssl}\colon C^{2,r}_0(\ov{\Omega})\to C^r(\ov{\Omega})$ is an isomorphism.

\medskip

Indeed, the injectivity follows by assumption and we focus on the surjectivity. Letting~$f\in C^\beta(\ov{\Omega})$, we consider
    \begin{align*}
       \mathcal D\varphi-\lm\rlq [\varphi]_{\ssl} =f,
    \end{align*}
    or equivalently
    \begin{align*}
        \varphi-T_{\ssl}(\varphi) = G*f.
    \end{align*}
    Using the projections in Proposition \ref{pro-proj} we get~$\varphi=\varphi_0+\varphi_1$, and~$G*f=(G*f)_0+ (G*f)_1$, with~$\varphi_0, (G*f)_0 \in \Eigen(T_{\ssl};0)$ and~$\varphi_1, (G*f)_1\in \mathcal{H}_1$. Observe also that~$T_{\ssl}(\varphi_0)=0$. Therefore,
    $$
		\varphi_0 =(G*f)_0,
		$$
		\begin{align}\label{eq:surjectivity of L}
        (I-T_{\ssl})\varphi_1= (G*f)_1.
    \end{align}
    Since~$1\notin\Spect(T_{\ssl})$, by the Fredholm alternative there exists a unique~$\varphi_1$ solving~\eqref{eq:surjectivity of L}.

\medskip

Let now $\lm<\lm_*(\om,p)$. We can thus use the above result  together with Lemma \ref{lemE1} to extend any such solution to a full $C^1$ curve of solutions with $\lm<\lm_*(\om,p)$. Suppose by contradiction there exists any other solution with $\lm<\lm_*(\om,p)$ not lying on the former curve. Reasoning as above we can extend also this solution to a full curve. Since for $\lm=0$ there is a unique solution, the two curves would intersect in $\lm=0$ yielding a bifurcation point, which is impossible.
\finedim

\

\subsection{\bf Qualitative behavior} \, \newline

We turn now to the monotonicity of the solutions, which is based on the properties of the spectral theory introduced in section \ref{sec:spectral}. Here, $\el$ is the energy of a solution, defined in \eqref{energy}.

\bpr\label{pr-enrg} Let $(\al_{\ssl},\psi_{\ssl})$ be a solution of \eqref{prob} such that $\sg_1(\all,\pl)> 0$. Then,
$$
\dfrac{d \el}{d\lm}>0.
$$
\epr

\proof
Since $0\notin\mbox{Spect}(L_{\ssl})$, locally around $(\al_{\ssl},\psi_{\ssl})$, the branch of solutions is a $C^1$ curve by Proposition \ref{pro-lin}. Let
$$
	\vl=\dfrac{d \pl}{d \lm}.
$$
By elliptic estimates, it is a $C^{2}_0(\ov{\om}\,)$ solution of
$$
\mathcal D \vl =\lm \rlq \vl + \rlq\pl +\rlq \frac{d\all}{d\lm}.
$$
The constraint in \eqref{prob} yields
$$
\frac{d\all}{d\lm}=-\lm  <\vl>_{\ssl}-<\pl>_{\ssl}.
$$
Therefore,
\beq\label{1b1}
\mathcal D \vl =\lm \rlq [\vl]_{\ssl} + \rlq[\pl]_{\ssl}.
\eeq
Now, by the definition of $E_{\lm}$, integrating by parts, and by using \eqref{prob} together with \eqref{1b1}, we have $E_{\ssl}=\frac12\ino |\nabla \pl|^2$, and in particular it holds,
$$
\frac{d}{d \lm}E_{\lm} = \ino \pl \mathcal D \vl= \ml \lm <\pl, [\vl]_{\ssl}>_{\ssl} +\ml   <\pl, [\pl]_{\ssl}>_{\ssl},
$$
where we recall that $\ml=\ino g'_{\lm}$. We first observe
$$
<\pl, [\pl]_{\ssl}>_{\ssl}=<[\pl]_{\ssl}, [\pl]_{\ssl}>_{\ssl}>0.
$$

\medskip

We consider here $\all<0$. The case $\all\geq0$ is easier and can be handled similarly. We use the decomposition in Proposition \ref{pro-proj} to write
\begin{align*}
    \pl=P_0(\pl)+\sum\limits_{j=1}^{+\ii}\beta_j\phi_j, \quad \vl=P_0(\vl)+\sum\limits_{j=1}^{+\ii}\gamma_j\phi_j,
\end{align*}
where
\begin{align*}
\beta_j   = \frac{\lm}{\mu_j}\ino \rlq[\phi_j]_{\ssl}[\pl]_{\ssl}, \quad
    \gamma_j    =\frac{\lm}{\mu_j} \ino \rlq [\phi_j]_{\ssl}[\vl]_{\ssl}.
\end{align*}

Multiplying \rife{1b1} by $\phi_j$, using \rife{lineq0} and integration by parts, we get
\begin{align*}\label{lamq31}
\sg_{j}\ino\rlq[\phi_j]_{\ssl}[\vl]_{\ssl} =
\ino\rlq [\phi_j]_{\ssl}[\pl]_{\ssl},
\end{align*}
from which we deduce~$\sg_j\gamma_j={\beta_j}$. Observe also that
~$\rlq[P_0(\pl)]_{\ssl}\equiv 0$ and similarly for $\vl$. Therefore,
\begin{align*}
    \ino\rlq[\vl]_{\ssl}[\pl]_{\ssl}
    &= \ino\rlq \bigr([P_0(\vl)]_{\ssl} + [P_1(\vl)]_{\ssl} \bigr) \bigr( [P_0(\pl)]_{\ssl} + [P_1(\pl)]_{\ssl}\bigr) \\
    &= \ino \rlq [P_1(\vl)]_{\ssl} [P_1(\pl)]_{\ssl} \\
    &= \ino \rlq \parenthesis{\sum_{j\geq 1} \gamma_j [\phi_j]_{\ssl}} \parenthesis{\sum_{k\geq 1} \beta_k [\phi_k]_{\ssl}}  \\
    &= \sum_{j\geq 1}\gamma_j \beta_j \ino \rlq [\phi_j]_{\ssl}^2.
\end{align*}
Recalling now that, by definition,
$$
  \ino \frac{\lm}{\mu_j} \rlq [\phi_j]_{\ssl}^2=1,
$$
we have
\begin{align*}	
    \sum_{j\geq 1}\gamma_j \beta_j \ino \rlq [\phi_j]_{\ssl}^2= \sum_{j\geq 1}  \gamma_j\beta_j \frac{\mu_j}{\lm}
    = \sum_{j\geq 1} \frac{\mu_j}{\lm} \sg_j \gamma_j^2.
\end{align*}
Finally, since by assumption~$\sg_1>0$,
$$
\sum_{j\geq 1} \frac{\mu_j}{\lm} \sg_j \gamma_j^2\geq \sg_1\sum_{j\geq 1} \frac{\mu_j}{\lm} \gamma_j^2=\sg_1\ino \rlq \parenthesis{\sum_{j\geq 1} \gamma_j[\phi_j]_{\ssl}}^2=\sg_1\ml\Abracket{[P_1(\vl)]_{\ssl}^2}_{\ssl}>0.
$$
This concludes the proof.
\finedim

\

We next turn to the monotonicity of $\all$ and prove the following. We recall here that $\nu_{1,\lm}$ stands for the standard first eigenvalue of $L_\lm$.
\bpr\label{pr3.2.best}
Let $(\al_{\ssl},\psi_{\ssl})$ be a solution of \eqref{prob} such that $\nu_{1,\lm}\geq 0$. Then,
$$
\dfrac{d \all}{d \lm}<0.
$$
\epr
\proof
By Lemma \ref{lem:xi1s} we know $0\notin\mbox{Spect}(L_{\ssl})$ and then, locally around $(\al_{\ssl},\psi_{\ssl})$, the branch of solutions is a $C^1$ curve by Proposition \ref{pro-lin}. Let $\xil=\lm \pl$, which satisfies
\beq \label{u}
\begin{cases}
&\mathcal D \xil =\lm g(x,\all+\xil) \quad \mbox{in}\;\;\om, \\\\
&\displaystyle{\ino} g(x,\all+\xil)=1
\end{cases}
\eeq
and let
$$
\wl=\frac{d\xil}{d\lm}
$$
which is a $C^{2}_0(\ov{\om}\,)$ function such that
\beq \label{w}
\mathcal D \wl =\lm \rlq [\wl]_{\ssl}+\rl= \lm\rlq\wl-\lm \rlq<\wl>_{\ssl}+\rl  \quad \mbox{in}\;\;\om.
\eeq
By the constraint in \eqref{u} we have
$$
\frac{d\all}{d\lm}=-<\wl>_{\ssl}.
$$
Suppose by contradiction that $<\wl>_{\ssl}\leq0$ and let $w_-=|\wl|\chi_{\om_-}$, where $\om_-=\{x\in\om \,:\, \wl(x)<0\}$. By a standard use of the Sard Lemma, we can assume w.l.o.g. that $\om_-$ is smooth. Multiply \eqref{w} by $w_-$ and integrate by parts to get
$$
	\ino (\mathcal D \wl)w_- = \lm\ino \rlq\wl w_- - \lm<\wl>_{\ssl}\ino\rlq w_-+\ino\rl w_-,
$$
$$
	-\ino \sum_{i,j}a_{ij}D_i w_- D_j w_- = -\lm\ino \rlq w^2_- - \lm<\wl>_{\ssl}\ino\rlq w_-+\ino\rl w_-
$$
and hence
\begin{align*}
\lm<\wl>_{\ssl}\ino\rlq w_- &= \ino \sum_{i,j}a_{ij}D_i w_- D_j w_- -\lm\ino \rlq w^2_- +\ino\rl w_- \\
	&\geq \nu_{1,\lm}\ino \rlq w^2_- +\ino\rl w_->0,
\end{align*}
since, by assumption, $\rl\geq0$ and $\rlq>0$ in the support of $\rl$, observing also that $\rl w_-\equiv0$ is not possible. Therefore, we deduce $<\wl>{\ssl}>0$, a contradiction.
\finedim

\

\

We can now prove Theorem \ref{thm2}.\\
{\em Proof of Theorem \ref{thm2}}.
The spectral estimate follows by Proposition \ref{preigen.LE}. The uniqueness part is contained in Proposition \ref{pro-lin}. The monotonicity of $\el$ and $\all$ is derived in Propositions \ref{pr-enrg} and \eqref{pr3.2.best}, respectively. This concludes the proof.
\finedim

\

\section{\bf Positivity threshold} \label{sec:positive}

\medskip

We collect here the discussion about the positivity threshold \eqref{positive}, starting by proving Theorem \ref{thm1}, i.e.
$$
 \lm_+(\om,p)>\frac{A}{p}\Lambda(\om,2p).
$$
{\em Proof of Theorem \ref{thm1}}.
Suppose by contradiction that $\lm_+\leq \frac{A}{P_0}\Lambda(\om,2p)$.
Then, since we know that $\lm_*>\frac{A}{P_0}\Lambda(\om,2p)$, we can use Proposition \ref{pro-lin} and Lemma \ref{lemE1}
to pass to the limit $\lm\nearrow\lm_+$ and obtain a solution $(\al_+,\psi_+)$ of  $\eqref{prob}_{\lm_+}$ with $\al_+\leq0$. Then
$u_+=\lm_+\psi_+$ satisfies, in particular
$$
\graf{\mathcal D u_+ = \lm_+ g(\al_++u_+) \quad \mbox{in} \;\;\om\\ \\
u_+>0 \;\;\mbox{in}\;\;\om, \quad u_+=0 \;\; \mbox{on}\;\;\pa\om.
}
$$
where with a little abuse of notation $g(u)=g(x,u)$. With the same notation, consider  the linearized problem
\beq\label{phis}
\graf{\mathcal D \phi = \lm_+ g'(\al_++u_+)\phi  \quad \mbox{in} \;\;\om\\ \\
\phi=0 \quad \mbox{on}\;\;\pa\om.
}
\eeq
We can now use the assumption \eqref{g1} on $g$ to deduce that $u_+$ is a positive strict subsolution of \rife{phis}. Then, it is well-known that the first eigenvalue of \rife{phis} is negative, i.e. $\nu_{1,\lm_+}(u_+)<0$. This contradicts Proposition \ref{preigen.LE} and concludes the proof.
\finedim

\medskip

We now consider the improved positivity estimates for the case $N=2$ under the extra assumption \eqref{extra} on the $x$ variable, i.e.
\beq \label{extra2}
C_1\,\tilde g(z) \leq g(\cdot,z)\leq C_2\,\tilde g(z).
\eeq
We also recall that here and in the sequel $\overline{\Omega_+}$ denotes the support of $g_{\lm}$. We start with the following energy estimate, of independent interest.

\bpr\label{altype2} Let $N=2$ and suppose \eqref{extra2} holds true. Let $(\all, \pl)$ be a solution of \eqref{prob} with $\all<0$.
Then,
\beq\label{level1.2}
\int\limits_{\om_+} \sum_{i,j}a_{ij} D_i\pl D_j\pl \leq \frac{1}{A}\frac{C_2}{C_1} \left(\frac{p+1}{8\pi}\right).
\eeq
\epr
\proof

Letting $\thl=\|\pl\|_{\ii}$ we introduce
$$
\om(t)=\{x\in\om\,:\,\pl>t\}, \quad \Gamma(t)=\{x\in\om\,:\,\pl=t\}, \quad t\in [0, \thl],
$$
and
$$
m(t)=\int\limits_{\om(t)}g(x,\al+\lm\pl),\qquad \mu(t)=|\om(t)|,
\qquad {e}(t)=\int\limits_{\om(t)}\sum_{i,j}a_{ij} D_i\pl D_j\pl.
$$
Here, $|\om(t)|$ stands for the area of $\om(t)$. The regularity of $m(t)$, $e(t)$ can be obtained by standard observations and we just refer to Lemma 2.2 in \cite{BJW2} for more details.
Observe that clearly
\beq\label{muem1}
m(0)=1, \quad \mu(0)=1, \quad e(0)=\int\limits_{\om}\sum_{i,j}a_{ij} D_i\pl D_j\pl,
\eeq
and
\beq\label{muem2}
m(\thl)=0, \quad \mu(\thl)=0, \quad e(\thl)=0.
\eeq
By the co-area formula and using assumption \eqref{extra2}
\begin{align*}
-m^{'}(t) &=\int\limits_{\Gamma(t)}\frac{g(x,\al+\lm\pl)}{|\nabla \pl |} \geq C_1\int\limits_{\Gamma(t)}\frac{\tilde g(\al+\lm\pl)}{|\nabla \pl |} \\
	&= C_1\,\tilde g(\al+\lm t)\int\limits_{\Gamma(t)}\frac{1}{|\nabla \pl |} = C_1\,\tilde g(\al+\lm t) (-\mu^{'}(t)),
\end{align*}
that is
\beq\label{diffm2}
-m^{'}(t)\geq C_1\,\tilde g(\al+\lm t) (-\mu^{'}(t)),
\eeq
for a.a. $t\in (t_0,\thl)$. Similarly,
\beq\label{diffm2+}
-m^{'}(t)\leq C_2\,\tilde g(\al+\lm t) (-\mu^{'}(t)).
\eeq
On the other hand, $m^{'}(t)= 0$ in $(0,t_0)$. Moreover,
\beq\label{diffem2}
m(t)=\int\limits_{\om(t)}\mathcal D \pl=\int\limits_{\Gamma(t)}\frac{\sum_{i,j}a_{ij} D_i\pl D_j\pl}{|\nabla \pl|}=-e^{'}(t),
\eeq
for a.a. $t\in (t_0, \thl)$ and $m(t)=m(0)=1$ for $t\in [0,t_0]$. We now use again \eqref{extra2} and apply Schwarz inequality to get
\begin{align}\nonumber
 -m^{'}(t)m(t)
 =&
\int\limits_{\Gamma(t)}\frac{g(x,\al+\lm\pl)}{|\nabla \pl |}\int\limits_{\Gamma(t)}\frac{\sum_{i,j}a_{ij} D_i\pl D_j\pl}{|\nabla \pl|} \\ \nonumber
\geq & AC_1\,\tilde g(\al+\lm t) \int\limits_{\Gamma(t)}\frac{1}{|\nabla \pl |}\int\limits_{\Gamma(t)}|\nabla \pl| \\ \nonumber
\geq & AC_1\,\tilde g(\al+\lm t)|{\Gamma(t)}|^2 \geq AC_1\,\tilde g(\al+\lm t)
4\pi \mu(t), \nonumber
\end{align}
for a.a. $t\in (t_0, \thl)$, where $|{\Gamma(t)}|$ is the length of $\Gamma(t)$ for which we use the isoperimetric inequality.
Therefore,
\beq\label{diff12}
\frac{(m^2(t))^{'}}{8\pi }+
AC_1\,\tilde g(\al+\lm t)\mu(t)\leq 0,
\eeq
for a.a. $t\in (t_0, \thl)$. Recalling the growth assumption \eqref{g2} on $g$, it is not difficult to check that
$$
\tilde g(\al+\lm t)\mu(t)\geq\frac{1}{\lm(p+1)}
\biggr(\tilde g(\al+\lm t)(\al+\lm t)\mu(t)\biggr)^{'}
-\frac{1}{\lm (p+1)}\tilde g(\al+\lm t)(\al+\lm t)\mu^{'}(t),
$$
$\mbox{for a.a. } t\in (t_0, \thl)$. Inserting the latte estimate into \rife{diff12} and using also \rife{diffm2+} we deduce
$$
\left(\frac{m^2(t)}{8 \pi}+
AC_1\frac{\tilde g(\al+\lm t)(\al+\lm t)}{\lm (p+1)}\mu(t) \right)^{'}
-A\frac{C_1}{C_2}\frac{1}{\lm (p+1)}\left( \all+\lm{t}\right)m^{'}(t)\leq 0,
$$
$\mbox{ for a.a. } t\in (t_0, \thl)$.
This yields
\beq\label{diff22}
-\frac{m^2(t)}{8\pi}-AC_1\frac{\tilde g(\al+\lm t)(\al+\lm t)}{\lm (p+1)}\mu(t)
-A\frac{C_1}{C_2}\frac{1}{\lm (p+1)}\int\limits_{t}^{\thl}ds \,\left( \all+\lm {s}\right)
m^{'}(s)\leq 0,
\eeq
for all $t\in (t_0,\thl)$. Recall now \rife{diffem2} and that $-\lm t_0=\all$. One has
\begin{align}\nonumber
 \int\limits_{t_0}^{\thl}ds \,\left( \all+\lm {s}\right)m^{'}(s)
 =&-\all m(t_0)+
\lm\int\limits_{t_0}^{\thl}ds\,s\, m^{'}(s)\\ \nonumber
=&-\all m(t_0)-\lm t_0m(t_0)-
\lm\int\limits_{t_0}^{\thl}ds m(s)\\ \nonumber
=& -\lm\int\limits_{t_0}^{\thl}ds m(s)=\lm\int\limits_{t_0}^{\thl}ds e^{'}(s)=-\lm e(t_0). \nonumber
\end{align}
Thus, for $t\to t_0^+$, \rife{diff22} reads
$$
-\frac{m^2(t_0)}{8\pi}+A\frac{C_1}{C_2}\frac{e(t_0)}{p+1}\leq 0,
$$
i.e.,
$$
e(t_0)\leq \frac{1}{A}\frac{C_2}{C_1}\frac{p+1}{8\pi}m^2(t_0)=\frac{1}{A}\frac{C_2}{C_1}\frac{p+1}{8\pi},
$$
which proves \rife{level1.2}.
\finedim

\

We derive now the improved positivity threshold.

\medskip

{\em Proof of Theorem \ref{thm3}}.
Suppose $\all<0$. We will prove
$$
\lm\geq A\left(\frac{C_1}{C_2}\right)^{\frac{p-1}{2p}}\left(\dfrac{8\pi}{p+1}\right)^{\frac{p-1}{2p}} \Lambda^{\frac{p+1}{2p}}(\om,p+1).
$$
Let us denote
$R_{p+1}(w)=\dfrac{\ino |\nabla w|^2}{\left(\ino |w|^{p+1}\right)^{\frac{2}{p+1}}}$, so that
$$
\Lambda(\om,p+1)=\inf\limits_{w\in H^{1}_0(\om)\setminus\{0\}}R_{p+1}(w).
$$
We will use the different expressions of the energy:
\begin{align*}
\el=&\frac12\ino \sum_{i,j}a_{ij} D_i\pl D_j\pl=\frac{1}{2\lm^2}\ino \sum_{i,j}a_{ij} D_i\xil D_j\xil \\
=&\frac{1}{2\lm}\ino g(\all+\xil)\xil=\frac{1}{2\lm}\int_{\om_+} g(\all+\xil)\xil.
\end{align*}
Hence we have
\begin{align}
 \lm=&\dfrac{\ino \sum_{i,j}a_{ij} D_i\xil D_j\xil}{\int_{\om_+}g(\all+\xil)\xil} \nonumber \\
 =& \dfrac{\int_{\om} \sum_{i,j}a_{ij} D_i\xil D_j\xil}{\int_{\om_+} \sum_{i,j}a_{ij} D_i\xil D_j\xil}
 \dfrac{\int_{\om_+} \sum_{i,j}a_{ij} D_i\xil D_j\xil}{\left(\int_{\om_+} \xil^{p+1}\right)^{\frac{2}{p+1}}}
 \dfrac{\left(\int_{\om_+} \xil^{p+1}\right)^{\frac{2}{p+1}}}
{\left(\int_\om g(\all+\xil)\xil\right)^{\frac{2}{p+1}}}
\frac{1}{\left(\int_\om g(\all+\xil)\xil\right)^{\frac{p-1}{p+1}}} \nonumber\\ \nonumber
\\ \nonumber
\geq& \dfrac{\int_{\om} \sum_{i,j}a_{ij} D_i\xil D_j\xil}{\int_{\om_+} \sum_{i,j}a_{ij} D_i\xil D_j\xil}
A\Lambda(\om_+,p+1)
\left(\dfrac{\int_{\om_+} \xil^{p+1}}{\int_{\om}g(\all+\xil)\xil}\right)^{\frac{2}{p+1}}
\frac{1}{\left(\lm\ino \sum_{i,j}a_{ij} D_i\pl D_j\pl \right)^{\frac{p-1}{p+1}}} \\ \nonumber
\\ \nonumber
=&\frac{\Lambda(\om_+,p+1)}
{\left(\lm \int_{\om_+} \sum_{i,j}a_{ij} D_i\pl D_j\pl \right)^{\frac{p-1}{p+1}}}
\dfrac{\int_{\om} \sum_{i,j}a_{ij} D_i\pl D_j\pl }{\int_{\om_+}\sum_{i,j}a_{ij} D_i\pl D_j\pl }
\left(\dfrac{\int_{\om_+} \xil^{p+1}}{\int_\om g(\all+\xil)\xil}\right)^{\frac{2}{p+1}}
\times \\ \nonumber
&\times\frac{\left( \int_{\om_+} \sum_{i,j}a_{ij} D_i\pl D_j\pl \right)^{\frac{p-1}{p+1}}}
{\left(\int_{\om} \sum_{i,j}a_{ij} D_i\pl D_j\pl\right)^{\frac{p-1}{p+1}}} \nonumber \\
& \\ \nonumber
\geq & \frac{\Lambda(\om_+,p+1)}{\left(\lm \int_{\om_+} |\nabla \pl|^2\right)^{\frac{p-1}{p+1}}}\left(\dfrac{\int_{\om_+} \xil^{p+1}}{\int_\om g(\all+\xil)\xil}\right)^{\frac{2}{p+1}} \\ \nonumber
& \\ \nonumber
 \geq &\frac{\Lambda(\om_+,p+1)}{\left(\lm \int_{\om_+} |\nabla \pl|^2\right)^{\frac{p-1}{p+1}}}\,, \nonumber
\end{align}
where in the last estimate we are using $\all<0$ and the assumption \eqref{g2} on $g$. Finally, by the energy estimate in
Proposition \ref{altype2}, we get
$$
\lm\geq A\left(\frac{C_1}{C_2}\right)^{\frac{p-1}{2p}}\left(\dfrac{8\pi}{p+1}\right)^{\frac{p-1}{2p}} \Lambda^{\frac{p+1}{2p}}(\om_+,p+1),
$$
and the thesis follows by the monotonicity of $\Lambda(\om,p)$ w.r.t. $\om$, see \cite{CRa1}.
\finedim

\

\section*{Appendix}

\medskip

We briefly introduce here the standard formulation of Grad-Shafranov type equations \cite{Frei,Kad,Stac} and their relation with problem \eqref{prob}. We use a similar notation from the well known work \cite{BeBr}, where for fixed $I>0$, one seeks a solution $(\gamma,v)\in\R\times C^{2,r}(\ov \om)$ to
\beq \label{prob2}
\begin{cases}
\mathcal D v = \mathcal G(x,v) \quad \mbox{in } \om, \\ \\
-\bigints_{\partial\om} \dfrac{\partial v}{\partial\nu} = I, \\ \\
v=\gamma \quad \mbox{on } \partial\om,
\end{cases}
\eeq
where $\mathcal G$ shares analogous properties as the one stated in the introduction. Here, $v$ is related to the flux function of the
plasma which occupies the region $\om_+ = \{x\in\om \,:\, v(x)>0\}.$
The quantity $\mathcal G(x,v)$ is related to the electric density current and $I$ is related to the total electric current through $\om$.

\medskip

We consider in this paper a class of models of the form $\mathcal G(x,v)= Ig(x,v)$. It is straightforward to see that the relation to the problem \eqref{prob} is given by
$$
\begin{cases}
\psi= \dfrac{v-\gamma}{I}\,,  \\
\lm=I.
\end{cases}
$$
It seems that to draw some connection between \eqref{prob} and \eqref{prob2} for more general models one would need some homogeneity property on $\mathcal G(x,v)$, which appears to be too restrictive. We thus postpone this aspect to future studies.

\

\end{document}